\documentclass[lfeqn,12pt]{article}
\usepackage{amssymb,bbm,pifont,mathtools,mathabx,mathrsfs,
}
\usepackage{minitoc}
 
\usepackage[nottoc]{tocbibind}

\usepackage{comment}

\usepackage{multirow}
\usepackage[multiple]{footmisc}
\usepackage{xcolor}
\usepackage[colorlinks=true]{hyperref}
\usepackage{marginnote}

\newcommand{\Giu}{{\bigskip\noindent}}

\newcommand{\noi}{{\noindent}}

\newtheorem{theorem}{Theorem}
\newtheorem{definition}[theorem]{Definition}
\newtheorem{proposition}[theorem]{Proposition}
\newtheorem{lemma}[theorem]{Lemma}
\newtheorem{remark}[theorem]{Remark}
\newtheorem{conjecture}[theorem]{Conjecture}
\newtheorem{question}{Question}
\newtheorem{sublemma}[theorem]{Sublemma}
\newtheorem{corollary}[theorem]{Corollary}
\newtheorem{assumption}[theorem]{Assumption}
\newtheorem{notationalrem}[theorem]{Notational Remark}

\newtheorem{tools}[subsection]{$\negsp\negsp$}
\newcommand\asm[1]{ \begin{assumption}\label{#1} }
\newcommand\easm{ \end{assumption} }
\newcommand\dfn[1]{ \begin{definition}\label{#1} }
\newcommand\dfntwo[2]{ \begin{definition}[#2]\label{#1} }
\newcommand\edfn{ \end{definition} }
\newcommand\rem[1]{ \begin{remark}\label{#1} \small \rm}
\newcommand\remtwo[2]{ \begin{remark}[#2]\label{#1} \rm}
\newcommand\erem{ \end{remark} }
\newcommand\thm[1]{ \begin{theorem}\label{#1}}
\newcommand\thmtwo[2]{ \begin{theorem}[#2]\label{#1}}
\newcommand\ethm{ \end{theorem} }
\newcommand\pro[1]{ \begin{proposition}\label{#1}}       
\newcommand\protwo[2]{ \begin{proposition}[#2]\label{#1}}
\newcommand\epro{ \end{proposition} }
\newcommand\conj[1]{ \begin{conjecture}\label{#1}}       
\newcommand\conjtwo[2]{ \begin{conjecture}[#2]\label{#1}}
\newcommand\econj{ \end{conjecture} }
\newcommand\quest[1]{ \begin{question}\label{#1}}       
\newcommand\questtwo[2]{ \begin{question}[#2]\label{#1}}
\newcommand\equest{ \end{question} }
\newcommand\lem[1]{ \begin{lemma}\label{#1}}
\newcommand\lemtwo[2]{ \begin{lemma}[#2]\label{#1}}
\newcommand\elem{ \end{lemma} }
\newcommand\sublem[1]{ \begin{sublemma}\label{#1}}
\newcommand\sublemtwo[2]{ \begin{sublemma}[#2]\label{#1}}
\newcommand\esublem{ \end{sublemma} }
\newcommand\cor[1]{ \begin{corollary}\label{#1}}
\newcommand\cortwo[2]{ \begin{corollary}[#2]\label{#1}}
\newcommand\ecor{ \end{corollary} }
\newcommand\notrem[1]{ \begin{notationalrem}\label{#1} \sl}
\newcommand\enotrem{ \end{notationalrem} }

\newcommand\average[1]{{ \left\langle #1 \right\rangle}}

\newcommand\equ[1]{{\rm (\ref{#1})}}
\newcommand\beq[1]{ \begin{equation}\label{#1} }
\newcommand{\eeq}{ \end{equation} }

\newcommand\beqa[1]{ \begin{eqnarray} \label{#1}}
\newcommand{\eeqa}{ \end{eqnarray} }
\newcommand{\beqano}{ \begin{eqnarray*} }
\newcommand{\eeqano}{ \end{eqnarray*} }

\newcommand{\proof}{\par\medskip\noindent{\bf Proof\ }}
\newcommand{\prooff[1]}{\par\medskip\noindent{\bf Proof of Theorem~\ref{#1}\ }}

\newcommand\warning[1]{\textcolor{red}{{#1}}}
\newcommand\gwarning[1]{\textcolor{blue}{{#1}}}
\newcommand{\ie}{{\it i.e.\  }}

\newcommand{\resp}{{\it resp.\ }}
\newcommand{\sss}{{\it iff\ }}

\newcommand{\qed}{\hskip.5truecm
            \vrule width 1.7truemm height 3.5truemm depth 0.truemm
            \par\Giu}
\newcommand{\qedeq}{\hskip.5truecm
            \vrule width 1.7truemm height 3.5truemm depth 0.truemm}

\newcommand\ovl[1]{ \overline {#1} }

\newcommand\genf[1]{ \mathscr E_{#1} }
\newcommand\sppp{ { p^-} }
\newcommand\sppf{ { p^+} }
\newcommand\sppdp{ {\dot p^-} }
\newcommand\sppdf{ {\dot p^+} }
\newcommand\tp{ { t^-} }
\newcommand\tf{ {t^+} }

\newcommand\su[1]{ \frac{1}{ {#1}} }

\newcommand{\dpr}{ {\partial}   }

\newcommand\eqby[1]{\stackrel{\equ{#1}}{=}}
\newcommand\leby[1]{\stackrel{\equ{#1}}{\le}}

\renewcommand{\Im}{{\rm \, Im\,}}

\newcommand{\negsp}{\hspace{-.04truecm}}
\newcommand\ex{\, e}

\newcommand{\G}{ {\Gamma}   }

\newcommand{\vae }{ {\varepsilon}   }

\renewcommand{\l}{ {\lambda}   }

\renewcommand{\r}{ {\rho}   }

\renewcommand{\o}{ {\omega}   }
\renewcommand{\O}{ {\Omega}   }
\newcommand{\torus}{ {\mathbb{ T}}   }
\renewcommand{\natural}{ {\mathbb{ N}}   }
\newcommand{\real}{ {\mathbb{ R}}   }
\newcommand{\integer}{ {\mathbb{ Z}}   }
\newcommand{\complex}{ {\mathbb { C}}   }

\newcommand{\rn}{ {\real^2}   }

\font\teneufm=eufm10
\font\seveneufm=eufm7
\font\fiveeufm=eufm5
\newfam\eufmfam
\textfont\eufmfam=\teneufm
\scriptfont\eufmfam=\seveneufm
\scriptscriptfont\eufmfam=\fiveeufm

\newcommand\appA[1]{\section{#1}\label{app:A}
\renewcommand{\theequation}{A.\arabic{equation}}
           \setcounter{equation}{0}
\renewcommand{\thetheorem}{A.\arabic{theorem}}
           \setcounter{theorem}{0}
                  }

\newcommand{\wt}{\widetilde}

\newcommand{\id}{{\mathtt {id}}}

\renewcommand\subset{\subseteq}

\begin{document}


\date{\small \today}

\title{
{\bf  Non co-preservation of the $1/2$ \& $1/(2l+1)$--rational caustics along deformations of circles
}\\
}
\author{
V. Kaloshin\footnote{Vadim Kaloshin: University of Maryland, College Park, MD, USA \& 
Institute of Science and Technology Austria, Am Campus~1, 3400 Klosterneuburg, Austria. 
\emph{vadim.kaloshin@gmail.com}}\ , 
\
C. E. Koudjinan\footnote{Comlan Edmond Koudjinan: Institute of Science and Technology Austria (IST Austria), Am Campus~1, 3400 Klosterneuburg,
Austria. \emph{edmond.koudjinan@ist.ac.at}
} 
}
\maketitle
\begin{abstract}
For any given positive integer $l$, we prove that every plane deformation of a circle which preserves the $1/2$ 
 and $1/(2l+1)$--rational caustics is trivial \ie the deformation consists only of similarities (rescalings plus isometries).
\end{abstract}

\noi

\section{Notations} 
\begin{itemize}
\item[\tiny $\bullet$] $\nabla f(t)\coloneqq f(t)-f(t-2\pi/3)$
\item[\tiny $\bullet$]$p\coloneqq p(t)$, $\dot p\coloneqq \dot p(t)$, $\sppp\coloneqq p(\tp)$, $\sppdp\coloneqq \dot p(\tp)$, $\sppf\coloneqq p(\tf)$, $\sppdf\coloneqq \dot p(\tf)$
\item[\tiny $\bullet$] $\mathscr F_I(f)(t)\coloneqq \sum_{k\in I}f_{k} \ex^{ikt}$, the projection on the Fourier's modes $k$ in $I\subset\integer$. If $I=n\integer$ for some $n\in\integer$, we write $\mathscr F_n(f)=\mathscr F_{n\integer}(f)$. \ 
\item[\tiny $\bullet$] Denote by $C^w_\r(\torus,\complex)$, the set of 
analytic function on the strip
$$
\torus_\r\coloneqq \{z\in\complex\ : \ |\Im z|< \r\}/2\pi \integer\,,
$$
endowed with the sup--norm $\|f\|_\r\coloneqq \sup_{\torus_\r}|f|$.
\end{itemize}

\section{Introduction}
A billiard is a mathematical modeling of the dynamic of a confined massless particle without friction and reflecting elastically on the boundary (without friction): the particle moves along a straight line with constant speed till it hits the boundary, then reflects off with reflection angle equals to the angle of incidence and follows the reflected straight line. It was introduced by G.D. Birkhoff \cite{birkhoff2020periodic} in 1920. 

\noi
A key notion in billiard dynamic is that of caustic. 

\dfn{defcaust}
A caustic of the billiard dynamic in a domain $\O$ is a curve $\mathfrak{C}$ with the property that any billiard trajectory that is once tangent to $\mathfrak{C}$ stays tangent to $\mathfrak{C}$ after each reflection on the boundary.
\edfn
Mather \cite{mather1982glancing} proves the non--existence of caustics  if the curvature of the boundary vanishes at one point. Thus, as far as caustics are concerned, we can focus on billiards in strictly convex domains; such billiards with at least $C^3$--boundary 
will be called Birkhoff billiards.\footnote{Observe that if $\O$ is not convex, then the billiard map is not
continuous; in this article we will be interested only in strictly convex domains. Moreover, as pointed out by Halpern \cite{halpern1977strange}, if the boundary is
not at least $C^3$, then the flow might not be complete.}
However, a caustic, if it exists, need neither be convex nor differentiable.

\noi
Nevertheless, according to KAM Theory \cite{lazutkin1973existence,kovachev1990invariant}, a positive measure set of convex differentiable  caustics which accumulates on the boundary and on which the
motion is smoothly conjugate to a rigid rotation do exists for Birkhoff billiards provided the boundary of the domain is sufficiently smooth. in general, the billiard dynamic induces naturally an orientation preserving circle homeomorphism on each convex caustic, which in particular admits a rotation number, also called rotation number of the caustic. In particular, a caustic is called rational (\resp irrational) if its rotation number is rational (\resp irrational). In this work, we are mainly concerned with rational caustics.

\dfn{def1}
Given $m,n\in\natural$ with $m\ge 2$, we call a caustic $n/m$--{\bf rational} if all the corresponding tangential billiards trajectories are periodic with the rotation number $n/m$. We denote by $\G_{n/m}(\O)$ the collection of all the $n/m$--rational caustics of $\O$.
\edfn

\noi
unlike irrational caustics which tends to be robust under perturbation according to KAM Theory, rational caustics tends to be quite rigid and, therefore, break up under perturbation. All the rational caustics may be even destroyed as show by  Pinto-de-Carvalho and Ram{\'\i}rez-Ros \cite{pinto2013non} who 
proved that can perturb an elliptic billiard and destroy all its rational caustics.

\noi
In contrast, Kaloshin and Ke Zhang \cite{kaloshin2018density} proved that can perturb a Birkhoff billiard table and create a new $1/q$--rational caustic for sufficiently larges $q$, provide the boundary of the table is $C^r$ with $r>4$. However, nothing is known for smalls $q$. Moreover, it is not also known if one can always perturb a sufficiently smooth Birkhoff billiard table and creates, simultaneously, many rational caustics.

\noi
On the other side, a natural question is:

\quest{qtab0}
Can one perturb a sufficiently smooth Birkhoff billiard table and (co--)preserve many of its rational caustics?
\equest
The question is still widely open, even in the simplest case of co--preservation of two rational caustics. Actually, the following intriguing conjecture has been made by Tabachnikov over ten years ago:
\conjtwo{tabash}{S. Tabachnikov}
In a sufficiently small $C^r$ ($r=2,\cdots,\infty, w$) neighborhood of the circle there is no other billiard domain of constant width and preserving $1/3$--caustics. 
\econj
It has been proven by J. Zhang \cite{zhang2019coexistence} that in the class of $\integer_2$--symmetric analytic deformation of the circle  with certain
Fourier decaying rate, any deformation of the circle that co--preserves $1/2$ and $1/3$-rational caustics  is necessarily an isometric transformation.\\

\noi
In the present paper, we settle the analytic deformative case of Conjecture~\ref{tabash}. We prove that, for any given positive integer $l$, if an analytic deformation of the circle co--preserves $1/2$ and $1/(2l+1)$--rational caustics then this deformation is trivial \ie consists only of circles (see Theorem~\ref{teo2} below). 

\noi
A bounded convex planar domain may be parametrized in various way, amongst which we have the parametrization with support function.

\section{Support function and some facts}
\noi
Given a bounded convex planar domain $\O$ with $C^1$ boundary $\dpr\O$ such that the origin of the cartesian coordinates is in the interior of $\O$, we denote by $p_\O\colon [0,2\pi)\to [0,\infty)$ the support function of $\dpr\O$. Denoting by $(x(t),y(t))$ the cartesian coordinates of the point on $\dpr\O$ corresponding to $(t,p_\O(t))$, we have\footnote{We refer the reader to \cite{resnikoff2015curves} for more details.}
\beq{eqcartsup}
\left\{
\begin{aligned}
&x(t)=p_\O(t)\cos t-\dot p_\O(t)\sin t\\
&y(t)=p_\O(t)\sin t+\dot p_\O(t)\cos t,
\end{aligned}
\right.
\eeq
where $\dot p_\O$ denotes the derivative of $p_\O$.

\noi
Given a supporting function $p$, we associate the generating function, denoted by $\genf p$, of the billiard map in the corresponding domain and given by
\beq{genfspp}
\begin{aligned}
\genf p(t,\tf)&\coloneqq \sqrt{(x(t)-x(\tf)^2+(y(t)-y(\tf)^2)}
   \eqby{eqcartsup} \big( p^2+\dot p^2+\sppf^2+(\sppdf)^2 -2p\sppf\cos(t-\tf)-\\
   &\qquad-2p\sppdf\sin(t-\tf)+2\dot p\sppf\sin(t-\tf)-2\dot p\sppdf\cos(t-\tf)\big)^{1/2}
\end{aligned}
\eeq

\noi
We have the following nice characterization of convex domains with $1/2$--rational caustics:
\lem{lem1}
A bounded convex domain $\O$ with $C^0$ boundary posseses a $2$--periodic caustic 
\sss its support function $p$ is of constant width:
$$
p(t)=\frac\o{2}+\sum_{k\in\integer} p^{(k)}\ex^{i(2k+1)t},\qquad t\in\torus,
$$
where $\o$ is the average width of $\O$ and $\{p^{(k)}\}_{k\in\integer}\subset\complex$.
\elem

\noi
The following error--function is the basis of the Lagrangian alternative approach proposed by Moser and Levi.
\dfn{def0}
Given a bounded convex domain $\O$ with $C^1$ boundary and support function $p$, $u\in C^0(\torus)$ and $m\in\natural$, we set
\beq{eq3Caus0}
E^m(p,u)\coloneqq \dpr_1\genf{p}(u, u^+)+\dpr_2\genf{p}( u^-,u),
\eeq
where $ u^\pm(t)\coloneqq u(t\pm\frac{2\pi}{m})$. 
\noi
For the sake of simplicity, we shall write $E$ for $E^3$. 
\edfn



\noi
Following Moser--Levi\cite{levi2001lagrangian}, we have the characterization:
\lem{Mos3pecaustEq}
Given $m\in\natural\setminus \{1\}$, a bounded convex domain $\O\subset \rn$ whose support function $p\in C^1(\torus)$ admits  a $1/m$--periodic caustic \sss there is a homeomorphism $u\colon \torus \to\torus$ 
such that
\beq{eq3Caus}
E^m(p,u)= \dpr_1\genf{p}(u, u^+)+\dpr_2\genf{p}( u^-,u)= 0.
\eeq
\elem

\noi
Let $\O_0\coloneqq \{(x,y)\in\rn:\ x^2+y^2\le 1\}$ be the unit disc and consider the one--parameter family $\O_\vae $ of deformation of $\O_0$ such that
$$
p_{\O_\vae}(t)=1+\vae p_1(t)+O(\vae^2),\qquad \mbox{for some}\quad p_1\in C^1(\torus),
$$ 
and let $\mathscr D \coloneqq \{\O_\vae :\ \vae \ge 0\}.$

\rem{rem0}
For any $\l>0$, the generating function of the  disc $\l\O_0$ of radius $\l$  is 
$
\genf{\l}(t,\tf)=\l\cdot\sqrt{2(1-\cos(t-\tf))} 
$.
Thus, for any $m\in \natural$,
$$
E^m(\l,\id)=0,
$$
\ie $\l \O_0$ possesses a $1/m$--rational caustic. 
\erem


\noi 
The following extends Lemma~\ref{lem1} to all natural numbers and 
is contained in Ramirez--Ros\cite{ramirez2006break}. We provide in $\S \ref{ramros}$ an alternative proof.
\thm{teo01}  Let $m\in\natural$ with $m\ge 2$ and $\vae >0$. 
Assume  $\O_\vae$ admits a $1/m$--rational caustic. Then $p_{1,km}=0$ for all $k\in\integer\setminus\{0\}$. 

\noi
Consequently, 
the set of $\O\in \mathscr D$ having a $1/m$--rational caustics is a submanifold of $\mathscr D$ of infinite codimension. 
\ethm


\section{Main result}
Denote by $\mathscr C$ the set of deformations $\O_\vae$ of the unit disc $\O_\vae$ within the class of strictly convex plane domains,  whose support function 
$p_{\vae}\in C^3(\torus)$ and is 
of the form: 
$
p_\vae(t)=1+\vae\,p_1+O(\vae^2)
$
 with $p_1\in C^w_\r(\torus,\real)$, for some $\r>0$. 
Then, the following holds.
\thm{teo2}   
Let $l\in\natural$ and $\O_\vae\in\mathscr C$ 
be a deformation 
of the unit disc. 
Assume, there exists $\vae_0>0$ such that for any $\vae\in[0,\vae_0)$, $\O_\vae$ possesses a $1/2$--rational and a $1/(2l+1)$--rational caustic. Then, the deformation $\O_\vae$ is trivial: $\O_\vae$ is a disc for any $\vae\in[0,\vae_0)$.
\ethm
\section{Proof of Theorem~\ref{teo01}\label{ramros}}
The proof of Theorem~\ref{teo01} will be deduce from the following Lemma.

\lem{pu1}
Let $m>2$, 
 $p_1^*\coloneqq \sum_{k\in\integer}  p_{1,k}^*\ex^{ikt}\in C^2(\torus)$ and $u_1^*\coloneqq \sum_{k\in \integer}  u_{1,k}^*\ex^{ikt}\in L^2(\torus)$. 
 Then, 
\beq{E1eqg3}
E^m(1+\vae p_1^*,\id+\vae u_1^*)=O(\vae^2),
\eeq
\sss for any $k\in\integer\setminus m\integer$,
\beq{conclem5}
u_{1,k}^*= a_{m,k}\, p_{1,k}^*  \quad \mbox{and} \quad \mathscr F_{m\integer\setminus\{0\}}(p_{1}^*)=0,
\eeq

\noi
where 
\begin{align*}
a_{m,k}
&\coloneqq{i \left(k \cot ^2\left(\frac{\pi  k}{m}\right)-\cot \left(\frac{\pi }{m}\right) \cot \left(\frac{\pi  k}{m}\right)\right)}\,.
\end{align*}

\noi
\elem

\proof  By Lemma~\ref{lemapA1} (see $\S \ref{AppSubA1}$), $E^m(1+\vae p_1^*,\id+\vae u_1^*)=O(\vae^2)$ \sss
$$
 \sin \left(\frac{\pi }{m}\right) (\dot p_ 1^++\dot p_ 1^-+2  \dot p_ 1+ u_ 1^++ u_ 1^--2  u_1)+(p_ 1^--   p_ 1^+)\cos \left(\frac{\pi }{m}\right) =0
$$
\ie for any $k\in\integer$,

$$
-2 i \sin \left(\frac{\pi }{m}\right) \sin ^2\left(\frac{\pi  k}{m}\right) u_{1,k}= \left(2 k \sin \left(\frac{\pi }{m}\right) \cos ^2\left(\frac{\pi  k}{m}\right)-\cos \left(\frac{\pi }{m}\right) \sin \left(\frac{2 \pi  k}{m}\right)\right) p_{1,k}\,,
$$
which, in turn, is equivalent to \equ{conclem5} as, for all $k\in m\integer\setminus\{0\}$,
$$
2 k \sin \left(\frac{\pi }{m}\right) \cos ^2\left(\frac{\pi  k}{m}\right)-\cos \left(\frac{\pi }{m}\right) \sin \left(\frac{2 \pi  k}{m}\right)=2 k \sin \left(\frac{\pi }{m}\right)\not=0. \qedeq
$$ 


\section{Proof of Theorems~\ref{teo2}\label{profmain2}}
 We start setting up some notation. We shall denote
 $$
{\tiny\bullet}\  t^\pm\coloneqq t\pm\frac{2\pi}{3}\;,\quad p_n^\pm\coloneqq p_n(t^\pm)\;,\quad u_n^\pm\coloneqq u_n(t^\pm)\;,\quad P_N\coloneqq \sum_{n=0}^N \vae^n p_n\,\quad 
U_N\coloneqq \sum_{n=0}^N \vae^n u_n,
$$
where $p_0\coloneqq 1$ and $u_0\coloneqq \id$. 

$$
{\tiny\bullet} \
E^m\left(P_N,U_N\right)=\sum_{k=-N}^\infty  E_{N,k}^m\;\vae^{N+k}
\hspace{9.5cm}\ 
$$
where
$$
E_{N,k}^m \coloneqq \su{(N+k)!}\frac{d^{N+k}}{d\vae^{N+k}}E^m\left(P_N,U_N\right)\bigg|_{\vae=0}\;.
$$
For $m=3$, here and henceforth, we will drop the superscript $m$ and write $E$ for $E^3$.

\noi
The following Lemma will be needed.
 \lem{folk}
 Let $f\in C^\o_\r(\torus,\complex)$, for some $\r>0$. 
 If 
\beq{aspflk}
\sum_{k\in \integer} f_k\,\ovl{ f_{k-n}}=0,\quad \mbox{for all}\ n\in\integer\setminus\{0\},
\eeq
then $f\equiv f_0$.
\elem  
\proof
%
 Set $g(z)\coloneqq f(z)\ovl{f(z)}$ and 
 consider the usual scalar product on $L^2(\torus)$:
$$
\average{u,v}\coloneqq \int_\torus u\ovl v .
$$
Fix $0<\r'<\r$. Then, for all $k\in\integer$,
\beq{expodecay}
|f_k|\le \|f\|_\r\, \ex^{-\r |k|},
\eeq
so that, 
\beq{absconvf}
\sup_{z\in \torus_{\r'}}\sum_{k\in\integer} |f_k\ex^{ikz}|\le \sum_{k\in\integer} |f_k|\ex^{\r' |k|}\leby{expodecay} \|f\|_\r\sum_{k\in\integer} \ex^{(\r'-\r) |k|}<\infty.
\eeq
Thus, 
\beq{fourf}
f(z)=\sum_{k\in\integer} f_k\ex^{ikz}\,,\qquad \mbox{on}\qquad  \torus_{\r'},
\eeq
and, therefore,
\beq{fourg}
g(z)=\sum_{k\in\integer} g_k\ex^{ikz}\,,\qquad \mbox{on}\qquad  \torus_{\r'}.
\eeq
Moreover, 
for any $n\in\integer\setminus\{0\}$,
\begin{align*}
g_{n}
	&=\average{g,\ex^{int}}\\
	&=\int_\torus (\sum_{k\in \integer} f_k\,\ex^{ikt})(\sum_{m\in \integer} \ovl f_m\, \ex^{-imt})\ex^{-int}\, dt \\
	&=\sum_{k\in \integer} f_k\,\ovl{ f_{k-n}}\\
	&\eqby{aspflk} 0.
\end{align*}
Consequently, 
$g\overset{\equ{fourg}}\equiv 0$ on ${\torus_{\r'}}$ \ie $|f|^2\big|_{\torus_{\r'}}\equiv g_0$, and this holds for all $0<\r'<\r$. Thus, $|f|^2\equiv g_0$ on ${\torus_{\r}}$. 
Then, the open mapping theorem yields $f\equiv f_0$. 
\qed
 
 \prooff[teo2] \\
 {$\bullet$ \bf Case $n=1$: }We argue by contradiction. 
Let 
$\O\in \mathscr D_{2,3}$ with support function $p(t)=1+\vae\, p_1+\vae^2p_2+O(\vae^3), 
$
where $p_1\in C^w_{\tilde{\r}}(\torus,\real)$, $p_2\in C^3(\torus)$, for some $\tilde{\r}>0$
 and for $\vae$ close to $0$. Without loss of generality, we can assume that
 \beq{contrass0}
 p_{1,-1}=p_{1,1}=0,\qquad\mbox{and}\qquad p_1\not\equiv0.
 \eeq 
 Then, by Lemma~\ref{Mos3pecaustEq}, there exists $u(t)=t+\sum_{n=1}^\infty \vae^n u_n(t),
$ with $\{u_n\}_{n\ge 1}\subset C^0(\torus)$ such that
\beq{CruPt}
E(p,u)=0.
\eeq
Also, observe that, by Lemma~\ref{lem1}, we have\footnote{By making the normalization $\mathscr F_{0}(p_n)=0$, $n\ge 1$.}
\beq{No2Harmo}
\mathscr F_{2\integer}(p_n)=0, \qquad \forall\, n\ge 1.
\eeq
\noi
We have 
$
0=E(p,u)=E(1+\vae p_1,\id+\vae u_1)+O(\vae^2)
$, so that $E(1+\vae p_1,\id+\vae u_1)=O(\vae^2)$. Thus, by Lemma~\ref{pu1}, we have, for any $k\in\integer\setminus 3\integer$,
\beq{conclem5bis}
u_{1,k}= a_{3,k}\, p_{1,k}  \quad \mbox{and} \quad \mathscr F_{3\integer\setminus\{0\}}(p_{1})=0.
\eeq
Therefore, Lemma~\ref{lemapA1} yields 
\beq{EqE11pf.}
E(1+\vae p_1,\id+\vae u_1)=E_{1,1}\,\vae^2+O(\vae^3).
\eeq 
%
Now, using Lemma~\ref{recFormu}, we have 
\beq{FormRecEN+1biS}
E(P_{2},U_{2})\eqby{FormRecEN+1} E(P_{1},U_{1})+\wt E_{2,0}\,\vae^{2}+O(\vae^{3})\eqby{EqE11pf.}(E_{1,1}+\wt E_{2,0})\vae^{2}+O(\vae^{3}),
\eeq
where
\beq{sdg001bisO}
\begin{aligned}
&\wt E_{2,0}	=\su4(-p_{2}^++p_{2}^-)+\frac{\sqrt{3}}{4}(\dot p_{2}^++\dot p_{2}^-+2\dot p_{2})+\frac{\sqrt{3}}{4}(u_{2}^++u_{2}^--2u_{2})\,.
\end{aligned}
\eeq
Hence,
\beq{eqEP2U2p}
0=E(u,p)=E(P_{2},U_{2})+O(\vae^{3})\eqby{FormRecEN+1biS}(E_{1,1}+\wt E_{2,0})\vae^{2}+O(\vae^{3}),
\eeq
 which implies 
 \beq{Etild20}
 E_{1,1}+\wt E_{2,0}=0.
 \eeq
  Thus, $\mathscr F_6(E_{1,1})=-\mathscr F_6(\wt E_{2,0})\warning{\overset{\equ{sdg001bisO},\equ{No2Harmo}}=0}$. Then, specializing \equ{F152ZE1100Apa2mp1} to $m=1$, we obtain, for all $n\in\integer\setminus\{0\}$,
 
  \beq{systmstp01}
  \begin{aligned}
  & \sum_{k\in \integer}k (k-n) p_{1,6k+1}\, p_{1,6(n-k)-1}=0. 
  \end{aligned}
  \eeq
 
 \noi
Now, consider the auxiliary function  $f(z)\coloneqq \sum_{k\in\integer} f_k\ex^{ikz}$ with 
$
f_k\coloneqq k \,p_{1,6k+1}.$ Then 
$f\in C^\o_\r(\torus,\real)$, where $\r\coloneqq \tilde{\r}/2$. 
Moreover, as $\ovl{f_{k-n}}\coloneqq (k-n) \,\ovl{p_{1,6(k-n)+1}}=(k-n) \,p_{1,6(n-k)-1}$,  the last relation in \equ{systmstp01} then reads: $\sum_{k\in\integer}f_k\,\ovl{ f_{k-n}}=0,$ for all $n\in\integer\setminus\{0\}$. Therefore, Lemma~\ref{folk} yields  $f\equiv f_0=0$ \ie $p_{1,6k+1}=0$ for all $k\in\integer\setminus\{0\}$. But then, as $p_{1,-1}=p_{1,1}\eqby{contrass0}0$, we would get $p_1\equiv 0$, which contradicts \equ{contrass0}. 

\noi
{$\bullet$ \bf General case $n\in\natural$ : } The proof in the general case is completely identical, up to two minor adjustments. The first one is the Cohomological equation \equ{systmstp01} which,  according to \equ{F152ZE1100Apa2mp1}, is in general:
\beq{systmstp01gen}
-\frac{i}{16}\sum_{k\in \integer}\sum_{r=1}^{4m+1} \mathscr P^m_r(n,k)\, p_{1,2(2m+1)k+r}\, p_{1,2(2m+1)(n-k)-r}=0,\quad \forall\, n\in\integer\setminus\{0\}.
\eeq
But, each of the polynomials $\mathscr P^m_r(n,k)$ splits:
$$
\mathscr P^m_r(n,k)= -16i\,(c_{m,r}^*)^2(k-c_{m,r}^{**})(n-k+c_{m,r}^{**}),
$$
where 
$$c_{m,r}^*= \sqrt{(2 m+1)^3 \sin \left(\frac{\pi }{2 m+1}\right) \cot ^2\left(\frac{\pi  r}{2 m+1}\right)}\quad \mbox{and}\quad 
c_{m,r}^{**}\coloneqq \frac{\cot \left(\frac{\pi }{2 m+1}\right) \tan \left(\frac{\pi  r}{2 m+1}\right)-r}{4 m+2}.$$
 Hence, the auxiliary function $f(z)\coloneqq \sum_{k\in\integer} f_k\ex^{ikz}$ should be defined by:\footnote{This is the second adjustment.} $
f_k\coloneqq c_{m,r}^*(k-z_{m,r}) \,p_{1,6k+r}.$
\qed

\section*{Appendix}
\addcontentsline{toc}{section}{Appendices}
\setcounter{section}{0}
\renewcommand{\thesection}{\Alph{section}} 

\appA{Reccurent formula for $p_n$ and $u_n$ and Taylor's series expansion of $E^m(P_N,U_N)$}
 \subsection{Expansion of $E^m(P_1,U_1)$ \label{AppSubA1}}
Let 
 $p_1\coloneqq \sum_{k\in\integer}  p_{1,k}\ex^{ikt}\in C^2(\torus)$ and $u_1\coloneqq \sum_{k\in \integer}  u_{1,k}\ex^{ikt}\in C^0(\torus)$.
Set
 $P_1\coloneqq 1+\vae p_1$, $U_1\coloneqq \id+\vae u_1$ and $ U_1^\pm(t)\coloneqq U_1(t\pm\frac{2\pi}{m})$. 

\lem{lemapA1}
 Given any integer $m\ge 2$, we have
\beq{EaApA1}
\begin{aligned}
&E^m(1+\vae p_1,\id+\vae u_1)=E_{1,0}^m\,\vae+E_{1,1}^m\,\vae^2 +O(\vae^3),
\end{aligned}
\eeq
with
\newpage
\begin{align*}
E_{1,0}^m&\coloneqq \frac{1}{2  } \left(\sin \bigg(\frac{\pi }{m}\right) (\dot p_ 1^++\dot p_ 1^-+2  \dot p_ 1+ u_ 1^++ u_ 1^--2  u_1)+(p_ 1^--   p_ 1^+)\cos \left(\frac{\pi }{m}\right)\bigg)\,,\\
E_{1,1}^m&\coloneqq - \frac{\csc^2 \left(\frac{\pi }{m}\right)}{32  }\left(-5 \cos \left(\frac{\pi }{m}\right) (p_1^+)^2+\cos \left(\frac{3 \pi }{m}\right) (p_1^+)^2+6 \sin \left(\frac{\pi }{m}\right)  u_1   p_1^+ -2 \sin \left(\frac{3 \pi }{m}\right)  u_1   p_1^+-\right.\\
&  -6 \sin \left(\frac{\pi }{m}\right)  u_1^+   p_1^+ +2 \sin \left(\frac{3 \pi }{m}\right)  u_1^+   p_1^+ +10 \sin \left(\frac{\pi }{m}\right) \dot p_1  p_1^+ +2 \sin \left(\frac{3 \pi }{m}\right) \dot p_1  p_1^++\\
& +6 \sin \left(\frac{\pi }{m}\right)  \dot p_1^+    p_1^+ -2 \sin \left(\frac{3 \pi }{m}\right)  \dot p_1^+    p_1^+ +4 \cos \left(\frac{\pi }{m}\right) \ddot p_1  p_1^+ -4 \cos \left(\frac{3 \pi }{m}\right) \ddot p_1  p_1^+-\\
& -\left(\cos \left(\frac{3 \pi }{m}\right)-5 \cos \left(\frac{\pi }{m}\right)\right)  (p_1^-)^2 -\cos \left(\frac{\pi }{m}\right)  (u_1^+)^2 +\cos \left(\frac{3 \pi }{m}\right)  (u_1^+)^2+\\
& +\cos \left(\frac{\pi }{m}\right)  (u_1^-)^2 -\cos \left(\frac{3 \pi }{m}\right)  (u_1^-)^2 +\cos \left(\frac{\pi }{m}\right)  (\dot p_1^+)^2 -\cos \left(\frac{3 \pi }{m}\right)  (\dot p_1^+)^2-\\
& -\cos \left(\frac{\pi }{m}\right)  (\dot p_1^-)^2 +\cos \left(\frac{3 \pi }{m}\right)  (\dot p_1^-)^2 +2 \cos \left(\frac{\pi }{m}\right)  u_1   u_1^+ -2 \cos \left(\frac{3 \pi }{m}\right)  u_1   u_1^+- \\
& -2 \cos \left(\frac{\pi }{m}\right)  u_1   u_1^-+2 \cos \left(\frac{3 \pi }{m}\right)  u_1   u_1^--2 \cos \left(\frac{\pi }{m}\right)  u_1^+  \dot p_1+2 \cos \left(\frac{3 \pi }{m}\right)  u_1^+  \dot p_1+\\
& +2 \cos \left(\frac{\pi }{m}\right)  u_1^- \dot p_1-2 \cos \left(\frac{3 \pi }{m}\right)  u_1^- \dot p_1+2 \cos \left(\frac{\pi }{m}\right)  u_1   \dot p_1^+  -2 \cos \left(\frac{3 \pi }{m}\right)  u_1   \dot p_1^++\\
& +2 \cos \left(\frac{\pi }{m}\right)  u_1^+   \dot p_1^+  -2 \cos \left(\frac{3 \pi }{m}\right)  u_1^+   \dot p_1^+  -2 \cos \left(\frac{\pi }{m}\right) \dot p_1  \dot p_1^+  +2 \cos \left(\frac{3 \pi }{m}\right) \dot p_1  \dot p_1^+ +\\
& +2  p_1  \left(-\left(\cos \left(\frac{3 \pi }{m}\right)-5 \cos \left(\frac{\pi }{m}\right)\right)  p_1^+ +\left(\cos \left(\frac{3 \pi }{m}\right)-5 \cos \left(\frac{\pi }{m}\right)\right)  p_1^--\right.\\
&-2 \sin \left(\frac{\pi }{m}\right) \left(-4  u_1  \sin ^2\left(\frac{\pi }{m}\right)+2  u_1^- \sin ^2\left(\frac{\pi }{m}\right)-\cos \left(\frac{2 \pi }{m}\right)  u_1^+ + u_1^++\right.\\
&\left.\left. +2 \cos \left(\frac{2 \pi }{m}\right) \dot p_1+6 \dot p_1-\cos \left(\frac{2 \pi }{m}\right)  \dot p_1^+  + \dot p_1^+  -\cos \left(\frac{2 \pi }{m}\right)  \dot p_1^-  + \dot p_1^-  \right)\right)-\\
& -2 \cos \left(\frac{\pi }{m}\right)  u_1   \dot p_1^-  +2 \cos \left(\frac{3 \pi }{m}\right)  u_1   \dot p_1^-  -2 \cos \left(\frac{\pi }{m}\right)  u_1^-  \dot p_1^-  +2 \cos \left(\frac{3 \pi }{m}\right)  u_1^-  \dot p_1^- +\\
& +2 \cos \left(\frac{\pi }{m}\right) \dot p_1  \dot p_1^-  -2 \cos \left(\frac{3 \pi }{m}\right) \dot p_1  \dot p_1^-  -12 \sin \left(\frac{\pi }{m}\right)  u_1^+  \ddot p_1+4 \sin \left(\frac{3 \pi }{m}\right)  u_1^+  \ddot p_1-\\
&-12 \sin \left(\frac{\pi }{m}\right)  u_1^- \ddot p_1+4 \sin \left(\frac{3 \pi }{m}\right)  u_1^- \ddot p_1-24 \sin \left(\frac{\pi }{m}\right) \dot p_1 \ddot p_1+8 \sin \left(\frac{3 \pi }{m}\right) \dot p_1 \ddot p_1-\\
&-12 \sin \left(\frac{\pi }{m}\right)  \dot p_1^+   \ddot p_1+4 \sin \left(\frac{3 \pi }{m}\right)  \dot p_1^+   \ddot p_1-12 \sin \left(\frac{\pi }{m}\right)  \dot p_1^-   \ddot p_1+4 \sin \left(\frac{3 \pi }{m}\right)  \dot p_1^-   \ddot p_1+\\
&+4  p_1^- \sin \left(\frac{\pi }{m}\right) \left(2  u_1  \sin ^2\left(\frac{\pi }{m}\right)-2  u_1^- \sin ^2\left(\frac{\pi }{m}\right)+\cos \left(\frac{2 \pi }{m}\right) \dot p_1+3 \dot p_1-\cos \left(\frac{2 \pi }{m}\right)  \dot p_1^-  +\right.
\end{align*}
\begin{align*}
&\left. + \dot p_1^- -2 \sin \left(\frac{2 \pi }{m}\right) \ddot p_1\right)-12 \sin \left(\frac{\pi }{m}\right)  u_1^+  \ddot p_1^++4 \sin \left(\frac{3 \pi }{m}\right)  u_1^+  \ddot p_1^+-12 \sin \left(\frac{\pi }{m}\right)  u_1^- \ddot p_1^-+\\
&\left.+4 \sin \left(\frac{3 \pi }{m}\right)  u_1^- \ddot p_1^-\right).
\end{align*}

\noi
In particular,
if $\mathscr F_{(2m+1)\integer\setminus\{0\}}(p_{1})=0$ and $u_{1,k}= a_{2m+1,k}\, p_{1,k}$, for all $k\in \integer\setminus (2m+1)\integer$, then


\begin{align}
&{\mathscr F_{2(2m+1)}(E_{1,1}^{2m+1})=   \sum_{n\in \integer}n \ex^{i2(2m+1)nt} \sum_{k\in \integer}\sum_{r=1}^{4m+1} \mathscr P^m_r(n,k)\, p_{1,2(2m+1)k+r}\, p_{1,2(2m+1)(n-k)-r}}\,.\label{F152ZE1100Apa2mp1}
\end{align}
\begin{align*}
\mathscr P^m_r(n,k) &\coloneqq 
c_{m,r}\left(\left(-1+e^{\frac{2 i \pi  r}{2 m+1}}\right) \cos \left(\frac{\pi }{2 m+1}\right)-i \left(1+e^{\frac{2 i \pi  r}{2 m+1}}\right) \sin \left(\frac{\pi }{2 m+1}\right) (k (4 m+2)+r)\right)\\
 & \left(\left(1+e^{\frac{2 i \pi  r}{2 m+1}}\right) \sin \left(\frac{\pi }{2 m+1}\right) ((4 m+2) (n-k)-r)-i \left(-1+e^{\frac{2 i \pi  r}{2 m+1}}\right) \cos \left(\frac{\pi }{2 m+1}\right)\right) \,,\\
 c_{m,r}&\coloneqq  \frac{-4 (2 m+1)  \csc \left(\frac{\pi }{2 m+1}\right)}{\left(-1+e^{\frac{2 i \pi  r}{2 m+1}}\right)^2}\,.
\end{align*}
\elem
\proof For the sake of simplicity, we shall give the proof for $m=1$; the general case follows word--by--word the same lines.

\noi
$(i)$ Indeed,
\begin{align}
&\genf{{P_{1}}}^2(U_{1}, U_{1}^+)=2(1-\cos\frac{2\pi}{m})+2\vae \bigg((1-\cos\frac{2\pi}{m})(p_1+p_1^+)-\sin\frac{2\pi}{m}(u_1-u_1^+)-\nonumber\\
&\hspace{7cm}-\sin\frac{2\pi}{m}(\dot p_1-\dot p_1^+)\bigg)+O(\vae^2),\label{eqqw1}
\end{align}
Thus,
\begin{align}
&\genf{{P_{1}}}^{-1}(U_{1}, U_{1}^+)\eqby{eqqw1} \left(2\sin\frac{2\pi}{m}\right)^{-1} -\vae\,\left(2\sin\frac{2\pi}{m}\right)^{-3} \bigg((1-\cos\frac{2\pi}{m})(p_1+p_1^+)-\nonumber\\
&\hspace{4cm}-\sin\frac{2\pi}{m}(u_1-u_1^+)-\sin\frac{2\pi}{m}(\dot p_1-\dot p_1^+) \bigg)+O(\vae^2),\label{eqqw1divmi2}
\end{align}
and, substituting $t$ by $t-2\pi/m$ in \equ{eqqw1divmi2}, we obtain
\begin{align}
&\genf{{P_{1}}}^{-1}(U_{1}^-, U_{1})=\left(2\sin\frac{2\pi}{m}\right)^{-1} -\vae\,\left(2\sin\frac{2\pi}{m}\right)^{-3} \bigg((1-\cos\frac{2\pi}{m})(p_1^-+p_1)-\nonumber\\
&\hspace{5.5cm}-\sin\frac{2\pi}{m}(u_1^--u_1)-\sin\frac{2\pi}{m}(\dot p_1^--\dot p_1) \bigg) +O(\vae^2). \label{eqqw2divmi2}
\end{align}
Moreover,
\begin{align}
&2\genf{{P_{1}}}(U_{1}, U_{1}^+)\dpr_1\genf{{P_{1}}}(U_{1}, U_{1}^+)=-2\sin\frac{2\pi}{m}+2\vae \bigg(\dot p_1-\dot p_1^+\cos\frac{2\pi}{m} +(u_1-u_1^+)\cos\frac{2\pi}{m}-\nonumber\\
&\hspace{8cm}-(p_1^++p_1+\ddot p_1)\sin\frac{2\pi}{m}\bigg)+O(\vae^2) \label{eqqw3}
\end{align}
and
\begin{align}
&2\genf{{P_{1}}}(U_{1}^-, U_{1})\dpr_2\genf{{P_{1}}}(U_{1}^-, U_{1})= 2\sin\frac{2\pi}{m}+2\vae \bigg(\dot p_1-\dot p_1^-\cos\frac{2\pi}{m} -(u_1^--u_1)\cos\frac{2\pi}{m}+\nonumber\\
&\hspace{8cm}+(p_1^-+p_1+\ddot p_1)\sin\frac{2\pi}{m}\bigg)+O(\vae^2).  \label{eqqw4}
\end{align}
Therefore, writing $E^m(P_1,U_1)=\genf{{P_{1}}}^{-1}(U_{1}, U_{1}^+)\genf{{P_{1}}}(U_{1}, U_{1}^+)\dpr_1\genf{{P_{1}}}(U_{1}, U_{1}^+)+$\\ $+\genf{{P_{1}}}^{-1}(U_{1}^-, U_{1})\genf{{P_{1}}}(U_{1}^-, U_{1})\dpr_2\genf{{P_{1}}}(U_{1}^-, U_{1})$ and using \equ{eqqw1divmi2}--\equ{eqqw4}, we obtain the formula of the first order term $E_{1,0}^m$  in \equ{EaApA1}.

\noi
 Similarly, expanding \equ{eqqw1}--\equ{eqqw4} up to the second order, one gets the formula of the second order term $E_{1,1}^m$  in \equ{EaApA1}. 
%
\noi
Then, 
simple computations yields the formula \equ{F152ZE1100Apa2mp1}.
\qed

\subsection{Recurrent formula for $E^m(P_N,U_N)$ for $N\ge 2$\label{subsecA3}}
We adopt the same notations as in $\S  \ref{profmain2}$
\lem{recFormu}
Let $m\ge 2$, $N\ge 1$, $p_1,\cdots, p_{N+1}\in C^2(\torus)$, and $u_1,\cdots, u_{N+1}\in L^2(\torus)$. Then, we have
\beq{FormRecEN+1}
E^m(P_{N+1},U_{N+1})=E^m(P_{N},U_{N})+\wt E_{N+1,0}^m\,\vae^{N+1}
+O(\vae^{N+2}),
\eeq
with
\beq{sdg001}
\begin{aligned}
&\wt E_{N+1,0}^m	=\frac{1}{2} \bigg(\sin \left(\frac{\pi }{m}\right) \left(\dot p_ {N + 1}^+ +\dot p_ {N + 1}^-+2  \dot p_ {N + 1}+ u_ {N + 1}^++ u_ {N + 1}^--2  u_ {N + 1}\right)+\\
&\hspace{3cm}+\cos \left(\frac{\pi }{m}\right) \left( p_ {N + 1}^-- p_ {N + 1}^+\right)\bigg)\,.
\end{aligned}
\eeq
\elem

\proof For the sake of simplicity, we shall give the proof for $m=3$; the general case follows word--by--word the same lines.

\noi
 We have on one hand,
\begin{align}
&\genf{{P_{N+1}}}^2(U_{N+1}, U_{N+1}^+)= P_{N+1}^2(U_{N+1})+\dot P_{N+1}^2(U_{N+1})+P_{N+1}^2(U_{N+1}^+)+\dot P_{N+1}^2(U_{N+1}^+)-\nonumber\\
	&\quad-2\bigg(P_{N+1}(U_{N+1})P_{N+1}(U_{N+1}^+)+\dot P_{N+1}(U_{N+1})\dot P_{N+1}(U_{N+1}^+) \bigg)\cos(U_{N}-U_{N}^+)+\nonumber\\
	&\quad+2\bigg(\dot P_{N+1}(U_{N+1})P_{N+1}(U_{N+1}^+)- P_{N+1}(U_{N+1})\dot P_{N+1}(U_{N+1}^+) \bigg)\sin(U_{N}-U_{N}^+)=\nonumber\\
	&\quad=\genf{{P_{N}}}^2(U_{N}, U_{N}^+)+\vae^{N+1}\bigg(3(p_{N+1}+p_{N+1}^+)-\sqrt{3}(u_{N+1}-u_{N+1}^+)-\sqrt{3}(\dot p_{N+1}-\dot p_{N+1}^+)\bigg)+\nonumber\\
	&\hspace{2cm}+ O(\vae^{N+2}). \label{eqmoinEnp1en}
\end{align}

\noi
Therefore,
\beq{genngennp1p}
\begin{aligned}
&\genf{{P_{N+1}}}^{-1}(U_{N+1}, U_{N+1}^+)= \genf{{P_{N}}}^{-1}(U_{N}, U_{N}^+)-\frac{1}{2}{\genf{{P_{N}}}^{-3}(U_{N}, U_{N}^+)}\vae^{N+1}\bigg(3(p_{N+1}+p_{N+1}^+)-\nonumber\\
	&\hspace{.5cm}-\sqrt{3}(u_{N+1}-u_{N+1}^+)-\sqrt{3}(\dot p_{N+1}-\dot p_{N+1}^+)\bigg)+ O(\vae^{N+2})\nonumber\\
	&=\genf{{P_{N}}}^{-1}(U_{N}, U_{N}^+)-\frac{\sqrt{3}}{18}\vae^{N+1}\bigg(3(p_{N+1}+p_{N+1}^+)-\sqrt{3}(u_{N+1}-u_{N+1}^+)-\sqrt{3}(\dot p_{N+1}-\dot p_{N+1}^+)\bigg)+\nonumber\\
	&\qquad+ O(\vae^{N+2})
\end{aligned}
\eeq
 and, replacing $t$ by $t-2\pi/3$ in the above formula, 
  we get
\beq{genngennp2p}
\begin{aligned}
&\genf{{P_{N+1}}}^{-1}(U_{N+1}^-, U_{N+1})= \genf{{P_{N}}}^{-1}(U_{N}^-, U_{N})-\frac{\sqrt{3}}{18}\vae^{N+1}\bigg(3(p_{N+1}^-+p_{N+1})-\sqrt{3}(u_{N+1}^--u_{N+1})-\nonumber\\
	&\hspace{4cm}-\sqrt{3}(\dot p_{N+1}^--\dot p_{N+1})\bigg)+ O(\vae^{N+2}).
\end{aligned}
\eeq
Furthermore,
\begin{align}
&2\genf{{P_{N+1}}}(U_{N+1}, U_{N+1}^+)\;\dpr_1 \genf{{P_{N+1}}}(U_{N+1}, U_{N+1}^+)= \big(P_{N+1}(U_{N+1})+\ddot P_{N+1}(U_{N+1})\big)\times\nonumber\\
	&\times\bigg(2\dot P_{N+1}(U_{N+1})-2\dot P_{N+1}(U_{N+1}^+)\cos(U_{N+1}- U_{N+1}^+)+2P_{N+1}(U_{N+1}^+)\sin(U_{N+1}- U_{N+1}^+) \bigg)\nonumber\\
	&\ =\genf{{P_{N}}}(U_{N}, U_{N}^+)\dpr_1 \genf{{P_{N}}}(U_{N}, U_{N}^+) +\vae^{N+1}\bigg(2\dot p_{N+1}+\dot p_{N+1}^+-\sqrt{3}\big( p_{N+1}  +p_{N+1}^++\ddot p_{N+1}\big)+\nonumber\\
	&\qquad+u_{N+1}^+-u_{N+1}\bigg)+O(\vae^{N+2})\,,\label{1stparEpt1}
\end{align}
and
\begin{align}
&2\genf{{P_{N+1}}}(U_{N+1}^-, U_{N+1})\;\dpr_2 \genf{{P_{N+1}}}(U_{N+1}^-, U_{N+1})= \big(P_{N+1}(U_{N+1})+\ddot P_{N+1}(U_{N+1})\big)\times\nonumber\\
	&\times \bigg(2\dot P_{N+1}(U_{N+1})-2\dot P_{N+1}(U_{N+1}^-)\cos(U_{N+1}^-- U_{N+1})-2P_{N+1}(U_{N+1}^-)\sin(U_{N+1}^-- U_{N+1}) \bigg)\nonumber\\
	&\ =\genf{{P_{N}}}(U_{N}^-, U_{N})\dpr_2 \genf{{P_{N}}}(U_{N}^-, U_{N}) +\vae^{N+1}\bigg(2\dot p_{N+1}+\dot p_{N+1}^--\sqrt{3}\big(- p_{N+1}  +p_{N+1}^--\ddot p_{N+1}\big)+\nonumber\\
	&\qquad+u_{N+1}^--u_{N+1}\bigg)+O(\vae^{N+2}).\label{1stparEpt2}
\end{align}

Thus, writing 
\begin{align*}
E(P_{N+1},U_{N+1})&=  \genf{{P_{N+1}}}^{-1}(U_{N+1}, U_{N+1}^+)\cdot \genf{{P_{N+1}}}(U_{N+1}, U_{N+1}^+)\dpr_1 \genf{{P_{N+1}}}(U_{N+1}, U_{N+1}^+)+\\
	&\ +\genf{{P_{N+1}}}^{-1}(U_{N+1}^-, U_{N+1})\cdot \genf{{P_{N+1}}}(U_{N+1}^-, U_{N+1})\dpr_2 \genf{{P_{N+1}}}(U_{N+1}^-, U_{N+1}),
\end{align*}
one gets the \equ{FormRecEN+1}.

\bibliographystyle{apa}
\bibliography{BibtexDatabase}
\end{document}